\documentclass[11pt]{amsart}
\usepackage{amsmath,amsfonts,graphicx}
\usepackage{amssymb}
\usepackage{graphicx, color}
\usepackage{epic,eepic}
\usepackage{color}

\renewcommand{\labelenumi}{(\theenumi)}
\newcommand{\QED}{Q.E.D.}

\makeatletter

\@addtoreset{equation}{section}
\newtheorem{dfn}{Definition}[section]
\newtheorem{thm}[dfn]{Theorem}
\newtheorem{prp}[dfn]{Proposition}
\newtheorem{rmk}[dfn]{Remark}
\newtheorem{lmm}[dfn]{Lemma}
\newtheorem{cor}[dfn]{Corollary}

\begin{document}

\title[Boundary of the horseshoe locus for the H\'enon family]
{Boundary of the horseshoe locus for the H\'enon family}

\author{Zin Arai}
\address{Chubu University Academy of Emerging Sciences, Kasugai, Aichi 487--8501, Japan.}
\email{zin@isc.chubu.ac.jp}

\author{Yutaka Ishii}
\address{Department of Mathematics, Kyushu University, Motooka, Fukuoka 819--0395, Japan.}
\email{yutaka@math.kyushu-u.ac.jp}

\author{Hiroki Takahasi}
\address{Department of Mathematics, Keio University, Hiyoshi, Yokohama 223--8522, Japan.}
\email{hiroki@math.keio.ac.jp}

\date{\textit{January 29, 2018.}}

\begin{abstract}
The purpose of this article is to investigate geometric properties of the parameter locus of the H\'enon family where the uniform hyperbolicity of a horseshoe breaks down. As an application, we obtain a variational characterization of equilibrium measures ``at temperature zero'' for the corresponding non-uniformly hyperbolic H\'enon maps. The method of the proof also yields that the boundary of the hyperbolic horseshoe locus in the parameter space consists of two monotone pieces, which confirms a conjecture in~\cite{AI}. The proofs of these results are based on the machinery developed in~\cite{AI} which employs the complexification of both the dynamical and the parameter spaces of the H\'enon family together with computer assistance.
\end{abstract}

\maketitle

\section{Statements of Results}

This is a follow-up article of the paper~\cite{AI}. In the current paper we continue our study on the dynamics of the H\'enon family:
\[f_{a, b} : (x, y)\longmapsto (x^2-a-by, x),\] 
where $(a, b)\in \mathbb{R}\times\mathbb{R}^{\times}$ is the parameter. The H\'enon map has been regarded as one of the most fundamental classes of non-uniformly hyperbolic systems. It is known that for certain choices of parameters, $f_{a, b}$ is topologically conjugate to a translation and hence has empty non-wandering set. It is also known that for certain choices of parameters, $f_{a, b}$ exhibits Smale's horseshoes on which the dynamics is uniformly hyperbolic~\cite{DN}. Therefore, the H\'enon family describes how a hyperbolic horseshoe dynamics is created through homoclinic or heteroclinic bifurcations.

The purpose of this article is to investigate geometric properties of the boundary of the the parameter locus for the H\'enon family $f_{a, b}$ where the uniform hyperbolicity of a horseshoe breaks down. As an application, we obtain a variational characterization of equilibrium measures ``at temperature zero'' for the corresponding non-uniformly hyperbolic H\'enon maps. The method of the proof also yields that the boundary consists of two monotone real-analytic pieces, which confirms a conjecture in~\cite{AI}. The proofs of these results are based on the machinery developed in~\cite{AI} which employs the complexification of both the dynamical and the parameter spaces of the H\'enon family together with computer assistance.

To state the results, let us first introduce two parameter loci. A H\'enon map $f_{a, b}$ is called a \textit{hyperbolic horseshoe on $\mathbb{R}^2$} if the restriction to its non-wandering set $f_{a, b} : \Omega(f_{a, b})\to \Omega(f_{a, b})$ is uniformly hyperbolic and topologically conjugate to the full shift of two symbols. 
Define
\[\mathcal{H}_{\mathbb{R}}\equiv \bigl\{(a, b) \in\mathbb{R}\times\mathbb{R}^{\times} : f_{a, b} \mbox{ is a hyperbolic horseshoe on } \mathbb{R}^2 \bigr\}\]
and call it the \textit{hyperbolic horseshoe locus} of $f_{a, b}$. It is also known~\cite{FM} that the topological entropy of $f_{a, b}$ satisfies $0\leq h_{\mathrm{top}}(f_{a, b})\leq \log 2$ for any $(a, b)\in\mathbb{R}\times\mathbb{R}^{\times}$. 
This leads us to define
\[\mathcal{M}_{\mathbb{R}}\equiv \bigl\{(a, b) \in\mathbb{R}\times\mathbb{R}^{\times} : f_{a, b} \mbox{ attains the maximal entropy}\, \log 2 \bigr\}\]
which we call the \textit{maximal entropy locus} of $f_{a, b}$. Note that by the continuity of the function $(a, b)\mapsto h_{\mathrm{top}}(f_{a, b})$ (see, e.g.~\cite{M}), we have $\overline{\mathcal{H}_{\mathbb{R}}}\subset\mathcal{M}_{\mathbb{R}}$ in $\mathbb{R}\times\mathbb{R}^{\times}$.

In the paper~\cite{AI} it has been shown that there exists a real analytic function $a_{\mathrm{tgc}} : \mathbb{R}^{\times}\to\mathbb{R}$ from the $b$-axis to the $a$-axis of the parameter space $\mathbb{R}\times\mathbb{R}^{\times}$ so that 
\begin{enumerate}
\renewcommand{\labelenumi}{(\roman{enumi})}
\item $(a, b)\in\mathcal{H}_{\mathbb{R}}$ iff $a>a_{\mathrm{tgc}}(b)$, 
\item $(a, b)\in\mathcal{M}_{\mathbb{R}}$ iff $a\geq a_{\mathrm{tgc}}(b)$,
\item when $a=a_{\mathrm{tgc}}(b)$, the map $f_{a, b}$ has exactly one orbit of homoclinic (resp. heteroclinic) tangencies of invariant manifolds of suitable fixed points when $b>0$ (resp. $b<0$),
\end{enumerate}
which extends a previous result of Bedford--Smillie~\cite{BS} for the case $|b|<0.06$ (see the Main Theorem in~\cite{AI}). 
Moreover, as a consequence of this result we conclude that both 
\[\mathcal{H}_{\mathbb{R}}^{\pm}\equiv \mathcal{H}_{\mathbb{R}} \cap \big\{(a, b) \in \mathbb{R}\times\mathbb{R}^{\times} : \pm b>0\big\} \quad \mbox{and} \quad \mathcal{M}_{\mathbb{R}}^{\pm}\equiv \mathcal{M}_{\mathbb{R}} \cap \big\{(a, b) \in \mathbb{R}\times\mathbb{R}^{\times} : \pm b>0 \big\},\] 
are connected and simply connected in $\big\{(a, b) \in \mathbb{R}\times\mathbb{R}^{\times} : \pm b>0\big\}$, that $\overline{\mathcal{H}_{\mathbb{R}}^{\pm}}=\mathcal{M}_{\mathbb{R}}^{\pm}$ and that $\partial \mathcal{H}_{\mathbb{R}}^{\pm}=\partial \mathcal{M}_{\mathbb{R}}^{\pm}$ (see the Main Corollary in~\cite{AI}).

The first main result of this article concerns a local geometric property of the loci boundary $\partial \mathcal{H}^+_{\mathbb{R}}=\partial \mathcal{M}^+_{\mathbb{R}}$ near $(a, b)=(2, 0)$. Namely,

\begin{thm} \label{thm:slope}
We have 
\[\frac{9}{8}<\lim_{b\to +0}\frac{da_{\mathrm{tgc}}}{db}(b)<\frac{23}{8}.\]
\end{thm}

Theorem \ref{thm:slope} is not just an estimate but can be applied to study certain ergodic properties of H\'enon maps $f=f_{a, b}$ at the first bifurcation parameters $(a, b)\in\partial\mathcal{H}^+_{\mathbb{R}}$. To state them, we here recall some terminologies. Let $f=f_{a, b}$ be a H\'enon map defined on $\mathbb{R}^2$ with $0<|b|<1$. Let $E^u_p$ be the one-dimensional subspace of $T_p\mathbb R^2$ satisfying 
\[\limsup_{n\to+\infty}\frac{1}{n}\log\|D_pf^{-n}|E_p^u\|<0.\]
Since $f^{-1}$ expands area, the one-dimensional subspace with this property is unique if it exists.
It was proved in \cite{ST} that $E^u_p$ is defined for all $p\in\Omega(f_{a, b})$.
Let $M(f)$ be the space of $f$-invariant Borel probability measures of a H\'enon map $f$ endowed with the weak topology. Let
\[\Lambda_{\nu}(f)\equiv \int\log \|D_pf|E^u_p\|d\nu(p),\]
be the unstable Lyapunov exponent of $\nu\in M(f)$.
Given $t\in\mathbb R$ an \textit{equilibrium measure} for the potential $-t\log \|D_pf|E_p^u\|$ is a measure which attains the supremum:
\[\sup\big\{h_{\nu}(f)-t\Lambda_\nu(f)\colon\nu\in M(f)\big\},\]
where $h_{\nu}(f)$ denotes the measure-theoretic entropy of $\nu$. Apart from the uniqueness, the existence of equilibrium measures for all $t>0$ was established in \cite{T}.

Since $t$ represents the inverse of temperature in statistical mechanics, it is natural to study the limit of the equilibrium measures as $t\to+\infty$.
An invariant measure $\mu\in M(f)$ is called a $(+)$-\textit{ground state} if there exists an increasing sequence $t_n\in\mathbb{R}$ with $t_n\to+\infty$ as $n\to\infty$ so that $\mu$ is obtained as the weak limit of equilibrium measures $\mu_n$ for the potential function $-t_n\log\|D_pf|E^u_p\|$. 
Since $\Omega(f_{a,b})$ is compact, ground states exist. 
To characterize ground states we introduce the following definition. A measure $\mu\in M(f)$ is called \textit{Lyapunov minimizing} if it satisfies $\Lambda(a, b)=\Lambda_{\mu}(f_{a, b})$, where
\[\Lambda(a, b)\equiv \inf_{\nu\in M(f_{a, b})}\Lambda_{\nu}(f_{a, b}).\] 
A measure $\mu\in M(f)$ is called \textit{entropy maximizing among the Lyapunov minimizing measures} if it attains the supremum of the metric entropy over all Lyapunov minimizing measures. 

Since the unstable Lyapunov exponent is not lower semi-continuous as a function of measures, the existence of Lyapunov minimizing measures is an issue. 
For instance, see~\cite{L} in which certain horseshoes with three symbols at the first bifurcation was shown to have no Lyapunov minimizing measure. 
On the other hand, a sufficient condition was introduced in~\cite{T} for the existence of Lyapunov minimizing measures of H\'enon-like maps at the first bifurcation. To state it, let $U_{\varepsilon}\subset \mathbb{R}^2$ is the $\varepsilon$-neighborhood of the Chebyshev parameter $(a, b)=(2, 0)$. 
When $\varepsilon>0$ is small, there is a saddle fixed point $Q=Q(a, b)$ of $f_{a, b}$ for $(a, b)\in U_{\varepsilon}$ obtained as the continuation of the fixed point $Q(2, 0)=(2, 2)$ of $f_{2, 0}$. 
Let $\lambda_Q(a, b)$ be the unstable eigenvalue for $D_Qf_{a, b}$. 
In Theorem A (a) of~\cite{T} the third-named author has proved that the non-degeneracy condition $\frac{1}{2}\log|\lambda_Q(a, b)| \ne \Lambda(a, b)$ implies the existence of Lyapunov minimizing measures. 
However, since $|\lambda_Q(a, b)|\to\log 4$ and $\Lambda(a, b)\to\log2$ as $(a, b)$ tends to $(2, 0)$ along $\partial\mathcal{H}^+_{\mathbb{R}}$, it was difficult to check the non-degeneracy condition for the H\'enon family by hand. 

As a consequence of Theorem \ref{thm:slope} we show that the non-degeneracy condition holds for the H\'enon family when $\varepsilon>0$ small. Namely,

\begin{thm} \label{thm:nondegenerate}
There exists $\delta>0$ so that the H\'enon map $f_{a, b}$ with $(a, b)\in\partial\mathcal{H}^+_{\mathbb{R}}\cap U_{\delta}$ satisfies $\frac{1}{2}\log|\lambda_Q(a, b)| > \Lambda(a, b)$.
\end{thm}

This theorem together with Theorem A (a) of~\cite{T} immediately yields the following variational characterization of invariant measures for H\'enon maps at the first bifurcation parameters.
 
\begin{cor} \label{cor:ergodic}
There exists $\delta>0$ so that any $(+)$-ground state of any H\'enon map $f_{a, b}$ with $(a, b)\in\partial\mathcal{H}^+_{\mathbb{R}}\cap U_{\delta}$ is Lyapunov minimizing, and entropy maximizing among the Lyapunov minimizing measures.
\end{cor}

The method of the proof of Theorem \ref{thm:slope} also yields the following global geometric property of the loci boundary $\partial \mathcal{H}_{\mathbb{R}}^{\pm}=\partial \mathcal{M}_{\mathbb{R}}^{\pm}$ which provides an affirmative answer to a conjecture in~\cite{AI} on the piecewise monotonicity of the function $a_{\mathrm{tgc}}$. 

\begin{thm} \label{thm:monotone}
The function $a_{\mathrm{tgc}} : \mathbb{R}^{\times}\to\mathbb{R}$ is strictly monotone decreasing on $\{b<0\}$ and strictly monotone increasing on $\{b>0\}$. Moreover, we have 
\[\lim_{b\to +0} \frac{da_{\mathrm{tgc}}}{d b}(b)\ne \lim_{b\to -0} \frac{da_{\mathrm{tgc}}}{d b}(b).\]
In particular, the Chebyshev parameter $(a, b)=(2, 0)$ is the unique corner of $\partial\mathcal{H}_{\mathbb{R}}\cup\{(2, 0)\}=\partial\mathcal{M}_{\mathbb{R}}\cup\{(2, 0)\}$ in the extended parameter space $(a, b)\in\mathbb{R}^2$.
\end{thm}

The proofs of Theorems \ref{thm:slope} and \ref{thm:monotone} heavily rely on the complexification of both the dynamical and the parameter spaces of the H\'enon family, together with some results in complex analytic geometry and complex dynamics. Therefore, Corollary \ref{cor:ergodic} indicates that a geometric property of a \textit{complex} parameter locus yields ergodic property of \textit{real} H\'enon maps.

Theorems \ref{thm:slope} and \ref{thm:nondegenerate} and Corollary \ref{cor:ergodic} are obtained by all the three authors of this paper. Theorem \ref{thm:monotone} is obtained by the first-named and the second-named authors of this paper.

\section{Proofs of Theorem \ref{thm:slope} and Theorem \ref{thm:monotone}}

In this section we prove Theorems \ref{thm:slope} and \ref{thm:monotone}. We first recall some terminologies from~\cite{AI}. 

Let $f_u$ and $f_v$ be two points in $\mathbb{C}^2$ which we call \textit{focuses}. 
Let $L_u$ (resp. $L_v$) be a complex line in $\mathbb{C}^2$ not containing $f_u$ (resp. $f_v$),  and let $L_u'$ (resp. $L_v'$) be the complex line through $f_u$ which is parallel to $L_u$ (resp. $L_v$). These define a pair of projections $\pi_u : \mathbb{C}^2\setminus L_u' \to L_u$ and $\pi_v : \mathbb{C}^2\setminus L_v' \to L_v$ which we call \textit{projective coordinates} in $\mathbb{C}^2$. A \textit{projective box} in $\mathbb{C}^2$ is a polydisk with respect to certain projective coordinates in $\mathbb{C}^2$. Let $\mathcal{B}=D_u\times_{\mathrm{pr}}D_v$ (resp. $\mathcal{B}'=D'_u\times_{\mathrm{pr}}D'_v$) be a projective box and let $(\pi_u, \pi_v)$ (resp. $(\pi'_u, \pi'_v)$) be the projective coordinates for $\mathcal{B}$ (resp. $\mathcal{B}'$). Let $f : \mathbb{C}^2 \to\mathbb{C}^2$ be a complex H\'enon map.

\begin{dfn}
We say that $f : \mathcal{B}\cap f^{-1}(\mathcal{B}')\to \mathcal{B}'$ satisfies the \textit{crossed mapping condition (CMC)} of degree $d$ if
\[\rho_f\equiv (\pi'_u \circ f, \pi_v \circ \iota) \, : \, \mathcal{B}\cap f^{-1}(\mathcal{B}') \longrightarrow D'_u\times D_v\]
is proper of degree $d$, where $\iota : \mathcal{B}\cap f^{-1}(\mathcal{B}')\to \mathcal{B}$ is the inclusion map. 
\end{dfn}

Let $\{\mathcal{B}_i\}_i$ be a family of projective boxes in $\mathbb{C}^2$. We set
\[\mathfrak{T}^+\equiv\bigl\{(0, 0), (0, 2), (0, 3), (1, 0), (2, 2), (2, 3), (3, 1) \bigr\}\]
and 
\[\mathfrak{T}^-\equiv\bigl\{(0, 0), (0, 2), (1, 0), (1, 2), (2, 4), (3, 4), (4, 1), (4, 3) \bigr\}.\]
Elements in $\mathfrak{T}^{\pm}$ are called \textit{admissible transitions}. A triple $(f, \{\mathcal{B}_i\}_i, \mathfrak{T}^{\pm})$ is said to satisfy the crossed mapping condition (CMC) if $f : \mathcal{B}_i\cap f^{-1}(\mathcal{B}_j)\to \mathcal{B}_j$ is a crossed mapping for any $(i, j)\in \mathfrak{T}^{\pm}$.

In the proofs of Theorems \ref{thm:monotone} and \ref{thm:slope} the following notion is crucial.
 
\begin{dfn}
An \textit{affine tin can} of center $(a_0, b_0)\in\mathbb{R}^2$, height $h>0$, radius $r>0$ and slope $s\in\mathbb{R}$ is
\[\mathcal{C}((a_0, b_0), h, r, s)\equiv\big\{(a, b)\in\mathbb{C}^2 : |b-b_0|\leq h,\, |a-a_0-s(b-b_0)|\leq r\big\}\]
and its vertical boundary is
\[\partial^v\mathcal{C}((a_0, b_0), h, r, s)\equiv \big\{(a, b)\in \mathbb{C}^2 : |b-b_0|\leq h,\, |a-a_0-s(b-b_0)|=r\big\}.\]
\end{dfn}

Let us first explain the proof of Theorem \ref{thm:slope} where we employ only one affine tin can, whereas the proof of Theorem \ref{thm:monotone} requires a family of affine tin cans. Let us write $\mathbb{D}_{\delta}\equiv \{|z|<\delta\}$. We use the following version of the Schwarz lemma.

\begin{lmm} \label{lmm:schwarz1}
Let $\varphi : \mathbb{D}_h\to\mathbb{D}_r$ be a holomorphic function with $\varphi(0)=0$. Then, 
\[\bigg|\frac{d\varphi}{dz}(0)\bigg|\leq \frac{r}{h}.\]
\end{lmm}

Set $\mathcal{E}\equiv \mathcal{C}((2, 0), 0.024, 0.021-\varepsilon, 2)$, where $\varepsilon=10^{-5}$. With computer assistance we show

\begin{prp}[Crossed Mapping] 
For $(a, b)\in\mathcal{E}$, there exists a family of boxes $\{\mathcal{B}^+_i\}_{i=0}^3$ so that $f_{a, b} : \mathcal{B}^+_i\cap f^{-1}_{a, b}(\mathcal{B}^+_j)\to \mathcal{B}^+_j$ is a crossed mapping for $(i, j)\in \mathfrak{T}^+$.
\label{prp:qt1}
\end{prp}

The proof of this proposition is identical to (iii) of Theorem 2.12 in~\cite{AI}, hence omitted.

Let $(a, b)\in\mathcal{E}$ and write $f=f_{a, b}$. Since $f : \mathcal{B}^+_0\cap f^{-1}(\mathcal{B}^+_0)\to \mathcal{B}^+_0$ is a crossed mapping of degree one by Proposition \ref{prp:qt1}, the local unstable manifold $V^u_{\mathrm{loc}}(Q)$ at $Q$ is a horizontal holomorphic disk of degree one in $\mathcal{B}^+_0$ and the local stable manifold $V^s_{\mathrm{loc}}(Q)$ at $Q$ is a vertical holomorphic disk of degree one in $\mathcal{B}^+_0$. As in~\cite{AI} this also enables us to define the notion of special pieces of $f$ as follows. Again thanks to Proposition \ref{prp:qt1}, 
\[V_{31\overline{0}}^s(a, b)^+\equiv \mathcal{B}^+_3\cap f^{-1}(\mathcal{B}^+_1\cap f^{-1}(V^s_{\mathrm{loc}}(Q)))\]
is a vertical holomorphic disk of degree one in $\mathcal{B}^+_3$, and
\[V_{\overline{0}23}^u(a, b)^+\equiv \mathcal{B}^+_3\cap f(\mathcal{B}^+_2\cap f(V^u_{\mathrm{loc}}(Q)))\]
is a horizontal holomorphic disk of degree two in $\mathcal{B}^+_3$. As described in Proposition 4.8 of~\cite{AI}, these pieces are responsible for the first bifurcation. This motivates to introduce

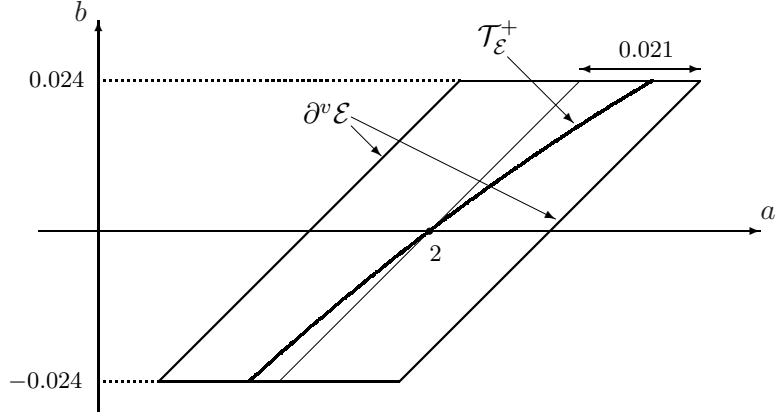
\begin{figure}
\setlength{\unitlength}{0.8mm}
\begin{picture}(120, 70)(0, 0)

\put(120, 32){$a$}
\put(6, 65){$b$}
\put(65, 30){\circle*{1}}
\put(73, 61){$\mathcal{T}^+_{\mathcal{E}}$}
\put(44, 48){$\partial^v\mathcal{E}$}

\footnotesize
\put(65, 25){$2$}
\put(96.5, 59){$0.021$}
\put(-5, 4){$-0.024$}
\put(-1.5, 54){$0.024$}

\thicklines
\put(20, 5){\line(1, 1){50}}
\put(60, 5){\line(1, 1){50}}
\qbezier(35,5)(68,35)(102,55)
\qbezier[10](10, 5)(15, 5)(20, 5)
\qbezier[60](10, 55)(40, 55)(70, 55)

\thinlines
\put(10, 0){\vector(0, 1){65}}
\put(0, 30){\vector(1, 0){120}}
\put(40, 5){\line(1, 1){50}}
\put(100, 57){\vector(1, 0){10}}
\put(100, 57){\vector(-1, 0){10}}
\put(78, 59){\vector(1, -1){11}}
\put(52, 47.5){\vector(1, -1){5}}
\put(52, 49){\vector(2, -1){34}}
\put(20, 5){\line(1, 0){40}}
\put(70, 55){\line(1, 0){40}}

\end{picture}
\caption{Complex tangency locus $\mathcal{T}_{\mathcal{E}}^+$ and the vertical boundary $\partial^v\mathcal{E}$.}
\label{FIG:tincan}
\end{figure}

\begin{dfn}
We define the \textit{complex tangency locus} in $\mathcal{E}$ as
\[\mathcal{T}_{\mathcal{E}}^+\equiv\bigl\{(a, b)\in\mathcal{E} : V_{31\overline{0}}^s(a, b)^+ \mbox{ and } V_{\overline{0}23}^u(a, b)^+ \mbox{ intersect tangentially}\bigr\}.\]
\end{dfn}

With computer assistance we show

\begin{prp}[Tin Can I] 
We have $\overline{\mathcal{T}_{\mathcal{E}}^+}\cap \partial^v \mathcal{E}=\emptyset$.
\label{prp:tincan1}
\end{prp}

See Figure \ref{FIG:tincan}. The proof of this proposition is identical to Theorem 5.4 in~\cite{AI}, hence omitted. Note, however, that the proof requires more accurate computation than (i) of Theorem 5.4 in~\cite{AI} because the width of $\mathcal{E}$ is much smaller compared with that of $\partial^v\mathcal{F}^+$ in~\cite{AI}.

\medskip

\noindent
\textit{Proof of Theorem \ref{thm:slope}.}
As in the proof of Proposition 5.9 in~\cite{AI}, Proposition \ref{prp:tincan1} together with the transversality of the quadratic family $p_a(x)=x^2-a$ at $a=2$ yields that $\mathcal{T}_{\mathcal{E}}^+$ is a complex submanifold of degree one in $\mathcal{E}$. Hence there exists a holomorphic function:
\[\kappa^+ : \big\{b\in\mathbb{C} : |b|<0.024\big\} \longrightarrow \mathbb{C}\]
whose graph coincides with $\mathcal{T}_{\mathcal{E}}^+$. Again by Proposition \ref{prp:tincan1},
\[\Psi : \big\{b\in\mathbb{C} : |b|<0.024\big\}\longrightarrow \big\{a\in\mathbb{C} : |a|<0.021-\varepsilon\big\}\]
given by $\Psi(b)\equiv \kappa^+(b)-2(b+1)$ is a well-defined holomorphic function. Lemma \ref{lmm:schwarz1} yields  
\[\biggl|\frac{d\kappa^+}{db}(0)-2 \biggr|=\bigg|\frac{d\Psi}{db}(0)\bigg|\leq\frac{0.021-\varepsilon}{0.024}<\frac{7}{8}.\]
Since $a_{\mathrm{tgc}}$ is defined as the real section of $\kappa^+$ for $0<b<1+\varepsilon$ in~\cite{AI}, we have
\[\frac{d\kappa^+}{db}(0)=\lim_{b\to +0}\frac{da_{\mathrm{tgc}}}{db}(b).\]
This finishes the proof of Theorem \ref{thm:slope}. \QED

\medskip

In the previous proof we employed only one affine tin can and the tangency locus $\mathcal{T}_{\mathcal{E}}^+$ goes through its center $(2, 0)$. The proof of Theorem \ref{thm:monotone} is based on a similar idea; we cover the tangency loci by a family of affine tin cans and apply the Schwarz lemma to conclude $\frac{da_{\mathrm{tgc}}}{db}(b)>0$ for $0< b\leq 1$ and $\frac{da_{\mathrm{tgc}}}{db}(b)<0$ for $-1\leq b< 0$. In this case, however, we need to construct an appropriate affine tin can for every $b$ and we do not know if the tangency locus goes through the centers of the affine tin cans. Therefore, a more ``flexible'' version of the Schwarz lemma is necessary.

\begin{lmm} 
Let $\varphi : \mathbb{D}_h\to\mathbb{D}_r$ be a holomorphic function. Then, 
\[\bigg| \frac{d\varphi}{dz}(z) \bigg|\leq \frac{h(r^2-|\varphi(z)|^2)}{r(h^2-|z|^2)}.\]
In particular, we have 
\[\bigg| \frac{d\varphi'}{dz}(z) \bigg| \leq \frac{4r}{3h}\]
for $|z|\leq \frac{h}{2}$.
\label{lmm:schwarz2}
\end{lmm}

For $0\leq n\leq 50$, let us set $b^{\pm}_n\equiv \pm 0.02\times n$. We then approximately compute the values $a_{\mathrm{tgc}}(b^{\pm}_n)$ by using the algorithm in~\cite{AM} and denote them by $a^{\pm}_n$ (see Tables \ref{tab:aprx1} and  \ref{tab:aprx2} at the end of this paper). Let $a_{\mathrm{aprx}} : \mathbb{R}\to\mathbb{R}$ be the piecewise affine interpolation of $a_{\mathrm{aprx}}(b^{\pm}_n)\equiv a^{\pm}_n$. It extends to $a_{\mathrm{aprx}} : \mathbb{C}\to\mathbb{R}$ by setting $a_{\mathrm{aprx}}(b)\equiv a_{\mathrm{aprx}}(\mathrm{Re}\,(b))$. We also approximately compute five values $a_{\mathrm{tgc}}(b^{\pm}_n-0.0002)$, $a_{\mathrm{tgc}}(b^{\pm}_n-0.0001)$, $a_{\mathrm{tgc}}(b^{\pm}_n)$, $a_{\mathrm{tgc}}(b^{\pm}_n+0.0001)$ and $a_{\mathrm{tgc}}(b^{\pm}_n+0.0002)$, and consider the degree four polynomial interpolation of these values. Let $s^{\pm}_n$ be the derivative of the degree four polynomial at $b=b^{\pm}_n$ (see Tables  \ref{tab:aprx1} and  \ref{tab:aprx2} again). Define
\[\mathcal{C}^{\pm}_n\equiv \mathcal{C}\bigg((a^{\pm}_n, b^{\pm}_n), h^{\pm}_n, \frac{3h^{\pm}_n}{2}, s^{\pm}_n\bigg),\]
where $h^+_n\equiv 0.01$ for $0\leq n\leq 50$ and
\[
h^-_n\equiv 
\begin{cases}
0.015 & \mbox{if } 19 \leq n \leq 21, \\
0.01 & \mbox{otherwise}.
\end{cases}
\]
Note that the union of the ``middle halves'' $\mathcal{C}\big((a^{\pm}_n, b^{\pm}_n), \frac{h^{\pm}_n}{2}, \frac{3h^{\pm}_n}{2}, s^{\pm}_n\big)$ of $\mathcal{C}^{\pm}_n$ forms a complex neighborhood of the graph of the function $a=a_{\mathrm{aprx}}(b)$ on $0<|b|\leq 1$.

With computer assistance we show

\begin{prp}[Quasi-Trichotomy] 
We have the following;
\begin{enumerate}
\renewcommand{\labelenumi}{(\roman{enumi})}
\item for any $1\leq n\leq 50$ and any $(a, b)\in\mathbb{R}^2$ satisfying $|b-b^{\pm}_n|\leq h^{\pm}_n$ and $a\leq a_{\mathrm{aprx}}(b)-\frac{3h^{\pm}_n}{2}$, we have $h_{\mathrm{top}}(f_{a, b}|_{\mathbb{R}^2})<\log 2$,
\item for any $1\leq n\leq 50$ and any $(a, b)\in\mathbb{R}^2$ satisfying $|b-b^{\pm}_n|\leq h^{\pm}_n$ and $a\geq a_{\mathrm{aprx}}(b)+\frac{3h^{\pm}_n}{2}$, $f_{a, b}$ is a hyperbolic horseshoe on $\mathbb{R}^2$,
\item for any $0\leq n\leq 50$ and any $(a, b)\in\mathcal{C}^{\pm}_n$, there exists a family of boxes $\{\mathcal{B}^{\pm}_i\}_i$ so that $f_{a, b} : \mathcal{B}^{\pm}_i\cap f^{-1}_{a, b}(\mathcal{B}^{\pm}_j)\to \mathcal{B}^{\pm}_j$ is a crossed mapping for $(i, j)\in \mathfrak{T}^{\pm}$.
\end{enumerate}
\label{prp:qt2}
\end{prp}

The proof of this proposition above is identical to Theorem 2.12 in~\cite{AI}, hence omitted. 

\begin{rmk}
In the proof of Theorem \ref{thm:monotone} we will see that it is not necessary to prove the claims (i) and (ii) of Proposition \ref{prp:qt2} for the case $n=0$. 
\end{rmk}

Thanks to (iii) of Proposition \ref{prp:qt2}, the special pieces $V_{31\overline{0}}^s(a, b)^+$ and $V_{\overline{0}23}^u(a, b)^+$ for $(a, b)\in \mathcal{C}^+_n$ and $V_{41\overline{0}}^s(a, b)^-$ and $V_{\overline{43}4124}^u(a, b)^-$ for $(a, b)\in \mathcal{C}^-_n$ can be defined for every $0\leq n\leq 50$.

\begin{dfn}\label{dfn:tangency}
We define the \textit{complex tangency loci} in $\mathcal{C}^{\pm}_n$ as
\[\mathcal{T}_n^+\equiv\bigl\{(a, b)\in\mathcal{C}^+_n : \mbox{$V_{31\overline{0}}^s(a, b)^+$ and $V_{\overline{0}23}^u(a, b)^+$ intersect tangentially}\bigr\}\]
and
\[\mathcal{T}_n^-\equiv\bigl\{(a, b)\in\mathcal{C}^-_n : \mbox{$V_{41\overline{0}}^s(a, b)^-$ and $V_{\overline{43}4124}^u(a, b)^-$ intersect tangentially}\bigr\}.\]
\end{dfn}

With computer assistance we show

\begin{prp}[Tin Can II]
We have $\overline{\mathcal{T}_n^{\pm}}\cap \partial^v\mathcal{C}^{\pm}_n=\emptyset$ for every $0\leq n\leq 50$.
\label{prp:tincan2}
\end{prp}

The proof of this proposition above is almost identical to Theorem 5.4 in~\cite{AI}, hence omitted. 

\begin{rmk}
The reason why we have to define $h^-_n$ discontinuously is the following. The complex subvariety $V_{\overline{43}4124}^u(a, b)^-$ in the definition of $\mathcal{T}_n^-$ (see Definition \ref{dfn:tangency}) consists of two irreducible components (see Proposition 5.11 in~\cite{AI}) which we are not able to ``distinguish'' a priori. One of them is responsible for the last tangency when $(a, b)$ is a real parameter, hence its tangency locus with $V_{41\overline{0}}^s(a, b)^-$ always belongs to the tin can $\mathcal{C}_n^-$. When $b<0$ is close to zero, the two components are close to each other and hence their tangency loci with $V_{41\overline{0}}^s(a, b)^-$ are both in $\mathcal{C}_n^-$. When $b<0$ is close to $-1$, the two components are so separated that one of the loci is outside $\mathcal{C}_n^-$. Therefore, if we keep $h^-_n$ to be a constant, we can never obtain a claim like Proposition \ref{prp:tincan2}.
\end{rmk}

Now let us finish the proof of our second main result.

\medskip

\noindent
\textit{Proof of Theorem \ref{thm:monotone}.}
As in the proof of Theorem \ref{thm:slope}, Proposition \ref{prp:tincan2} with Lemma \ref{lmm:schwarz2} imply 
\begin{equation} 
\bigg|\frac{d\kappa^{\pm}}{db}(b)-s_n^{\pm}\bigg| \leq \frac{4}{3}\cdot \frac{\frac{3h^{\pm}_n}{2}}{h^{\pm}_n}= 2
\label{eqn:slope}
\end{equation}
for $|b-b^{\pm}_n|\leq \frac{h^{\pm}_n}{2}$. This estimate for $n=0$ yields 
\begin{equation}
\frac{da_{\mathrm{tgc}}}{db}(b) \in [s^+_0-2, s^+_0+2] \quad \mbox{and} \quad \frac{da_{\mathrm{tgc}}}{db}(b) \in [s^-_0-2, s^-_0+2]
\label{eqn:initial}
\end{equation}
for $b \in \big(0, b^+_0+\frac{h^+_0}{2}\big]$ and for $b \in \big[b^-_0-\frac{h^-_0}{2}, 0\big)$ respectively. Note that here we did not use the claims (i) and (ii) of Proposition \ref{prp:qt2} for the case $n=0$. Although the function $a_{\mathrm{tgc}}$ is not defined at $b=0$, the estimates (\ref{eqn:initial}) hold even for $b=0$ in the sense that 
\[\lim_{b\to +0}\frac{da_{\mathrm{tgc}}}{db}(b)\in [0.000001, 4.000001] \quad \mbox{and} \quad \lim_{b\to -0}\frac{da_{\mathrm{tgc}}}{db}(b)\in [-4.246320, -0.246320]\]
(see Tables  \ref{tab:aprx1} and  \ref{tab:aprx2} again for the values of $s^{\pm}_0$). In particular, $\partial \mathcal{H}_{\mathbb{R}}\cup\{(2, 0)\}$ forms a corner at $(a, b)=(2, 0)$. 

Next, let us consider the case $1\leq n\leq 50$. By Proposition \ref{prp:tincan2}, $\mathcal{T}^{\pm}_n$ are complex subvarieties in $\mathcal{C}^{\pm}_n$ of some degrees. Since $\mathcal{T}^{\pm}_n$ are subsets of the complex tangency loci $\mathcal{T}^{\pm}$ introduced in~\cite{AI} and we know that $\mathcal{T}^{\pm}$ are degree one~\cite{AI}, the degrees of $\mathcal{T}^{\pm}_n$ are at most one. On the other hand, since the graph of the function $a=a_{\mathrm{tgc}}(b)$ have to go through $\mathcal{C}^{\pm}_n$ from its bottom boundary to the top boundary thanks to (i) and (ii) of Proposition \ref{prp:qt2}, the degrees of $\mathcal{T}^{\pm}_n$ are at least one. It follows that the degrees of $\mathcal{T}^{\pm}_n$ are all one. Then, we can proceed as in the case $n=0$. Namely, by Proposition \ref{prp:tincan2}, we obtain the estimate (\ref{eqn:slope}) for $|b-b^{\pm}_n|\leq \frac{h^{\pm}_n}{2}$ and for $1\leq n\leq 50$. This yields 
\begin{equation}
\frac{da_{\mathrm{tgc}}}{db}(b) \in [s^+_n-2, s^+_n+2] \quad \mbox{and} \quad \frac{da_{\mathrm{tgc}}}{db}(b) \in [s^-_n-2, s^-_n+2]
\label{eqn:estimate}
\end{equation}
for $b \in \big[b^+_n-\frac{h^+_n}{2}, b^+_n+\frac{h^+_n}{2}\big]$ and for $b \in \big[b^-_n-\frac{h^-_n}{2}, b^-_n+\frac{h^-_n}{2}\big]$ respectively. Since 
\[\bigg(0, b^+_0+\frac{h^+_0}{2}\bigg]\cup \bigcup_{n=1}^{50}\bigg[b^+_n-\frac{h^+_n}{2}, b^+_n+\frac{h^+_n}{2}\bigg] \supset (0, 1] \quad \mbox{and} \quad \bigg[b^-_0-\frac{h^-_0}{2}, 0\bigg)\cup \bigcup_{n=1}^{50}\bigg[b^-_n-\frac{h^-_n}{2}, b^-_n+\frac{h^-_n}{2}\bigg] \supset [-1, 0)\]
hold, the estimates (\ref{eqn:initial}) and (\ref{eqn:estimate}) yield 
\begin{equation}
\frac{da_{\mathrm{tgc}}}{db}(b) \in \bigcup_{n=0}^{50}[s^+_n-2, s^+_n+2]=[s^+_0-2, s^+_{50}+2]=[0.000001, 7.699311]
\label{eqn:derivative1}
\end{equation}
for all $0< b\leq 1$ and 
\begin{equation}
\frac{da_{\mathrm{tgc}}}{db}(b) \in \bigcup_{n=0}^{50}[s^-_n-2, s^-_n+2]=[s^-_{50}-2, s^-_0+2]=[-8.198261, -0.246320]
\label{eqn:derivative2}
\end{equation}
for all $-1 \leq b < 0$ (see Tables  \ref{tab:aprx1} and  \ref{tab:aprx2} again for the values of $s^{\pm}_0$ and $s^{\pm}_{50}$). 

Our final task is to extend the above monotonicity result of $a_{\mathrm{tgc}}$ for $0<|b|\leq 1$ to all $b\in\mathbb{R}^{\times}$. For this we use the fact that $f_{a, b}^{-1}$ is affine conjugate to $f_{a/b^2, 1/b}$. This motivates to consider the following transformation in the parameter space:
\[\Phi : \mathbb{R}\times\mathbb{R}^{\times} \ni (a, b) \longmapsto \bigg(\frac{a}{b^2}, \frac{1}{b}\bigg) \in \mathbb{R}\times\mathbb{R}^{\times}.\]
This transformation maps the parameter region $\{0<|b|\leq 1\}$ to $\{|b|\geq 1\}$. Note that the boundary of the hyperbolic horseshoe locus $\partial\mathcal{H}_{\mathbb{R}}^{\pm}$ is invariant under this transformation. The projectivization of the derivative $D_{(a, b)} \Phi : T_{(a, b)}\partial\mathcal{H}_{\mathbb{R}}^{\pm}\rightarrow T_{\Phi(a, b)}\partial\mathcal{H}_{\mathbb{R}}^{\pm}$ is computed as
\begin{equation}
\mathbb{P}(D_{(a, b)} \Phi) : \mathbb{P}(T_{(a, b)}\partial\mathcal{H}_{\mathbb{R}}^{\pm})\ni [t : 1]\longmapsto \bigg[\frac{2a}{b}-t : 1\bigg]\in\mathbb{P}(T_{\Phi(a, b)}\partial\mathcal{H}_{\mathbb{R}}^{\pm}).
\label{eqn:formula}
\end{equation}
First, by (i) and (ii) of Quasi-Trichotomy we see 
\begin{equation}
\frac{2a}{b}\in [2\times 5.699311, +\infty) \quad \mbox{and} \quad \frac{2a}{b}\in (-\infty, 2\times -6.198261]
\label{eqn:large1}
\end{equation}
for $(a, b)\in\partial\mathcal{H}_{\mathbb{R}}^+ \cap \{0 < b \leq 1\}$ and for $(a, b)\in\partial\mathcal{H}_{\mathbb{R}}^-\cap \{-1 \leq b < 0\}$ respectively. On the other hand, (\ref{eqn:derivative1}) and (\ref{eqn:derivative2}) yield
\begin{equation}
t \in [0.000001, 7.699311] \quad \mbox{and} \quad t \in [-8.198261, -0.246320]
\label{eqn:large2}
\end{equation}
for $(t, 1)\in T_{(a, b)}\partial\mathcal{H}_{\mathbb{R}}^+$ with $0 < b \leq 1$ and for $(t, 1)\in T_{(a, b)}\partial\mathcal{H}_{\mathbb{R}}^-$ with $-1 \leq b < 0$ respectively. By (\ref{eqn:large1}) and (\ref{eqn:large2}), we see $\frac{2a}{b}-t>0$ for $0 < b\leq 1$ and $\frac{2a}{b}-t<0$ for $-1\leq b < 0$. This demonstrates the monotonicity of $a_{\mathrm{tgc}}$ for $b\geq 1$ and for $b\leq -1$ by the formula (\ref{eqn:formula}). \QED

\section{Proofs of Theorem \ref{thm:nondegenerate} and Corollary \ref{cor:ergodic}}

In this section we deduce Theorem \ref{thm:nondegenerate} from Theorem \ref{thm:slope} and then prove Corollary \ref{cor:ergodic}. 

For $(a, b)\in U_{\varepsilon}$ let $P=P(a, b)$ be the fixed point of $f_{a, b}$ obtained as the continuation of the fixed point $P(2, 0)=(-1, -1)$ of $f_{2, 0}$. 
Let $\lambda_P(a, b)$ be the unstable eigenvalue of $D_Pf_{a, b}$. 

\medskip

\noindent
\textit{Proof of Theorem \ref{thm:nondegenerate}.}
Define
\[\Gamma(a, b)\equiv \frac{1}{2}\log|\lambda_Q(a, b)|-\log|\lambda_P(a, b)|.\]
An easy computation shows
\[\Gamma(a, b)= \frac{1}{2}\log \Big|x_Q+\sqrt{(1+b)x_Q+a-b}\Big|-\log \Big|x_P-\sqrt{(1+b)x_P+a-b}\Big|,\]
where
\[x_P=\frac{1+b-\sqrt{(1+b)^2+4a}}{2} \quad \textrm{and} \quad x_Q=\frac{1+b+\sqrt{(1+b)^2+4a}}{2}.\]
By the implicit function theorem, the zero locus $\big\{(a, b)\in U_{\delta} : \Gamma(a, b)=0\big\}$ forms a smooth curve $\gamma(b)=a$ through $(a, b)=(2, 0)$ for $\delta>0$ sufficiently small. Moreover, by using
\[\Gamma(a, 0)= \frac{1}{2}\log \Bigg|\frac{1+\sqrt{1+4a}}{2}+\sqrt{\frac{1+\sqrt{1+4a}}{2}+a}\Bigg|-\log \Bigg|\frac{1-\sqrt{1+4a}}{2}-\sqrt{\frac{1-\sqrt{1+4a}}{2}+a}\Bigg|,\]
one can check $\Gamma(2, 0)=0$ and $\frac{\partial\Gamma}{\partial a}(2, 0)=-\frac{1}{4}<0$, hence $\Gamma(a, 0)<0$ for $a>2$ close to $2$. 

On the other hand, we compute 
\[\frac{d\gamma}{db}(0)=-\frac{\frac{\partial\Gamma}{\partial b}(2, 0)}{\frac{\partial\Gamma}{\partial a}(2, 0)}=\frac{23}{8}.\]
Therefore, Theorem \ref{thm:slope} implies
\[\lim_{b\to +0}\frac{da_{\mathrm{tgc}}}{db}(b) < \frac{d\gamma}{db}(0).\]
This, together with $\Gamma(a, 0)<0$ for $a>2$ close to $2$, yields that there exists $\delta>0$ with $\Gamma(a, b)>0$ for $(a, b)\in \partial\mathcal{H}^+_{\mathbb{R}}\cap U_{\delta}$. Since we see
\[\frac{1}{2}\log|\lambda_Q(a, b)|-\Lambda(a, b) \geq \frac{1}{2}\log|\lambda_Q(a, b)|-\log|\lambda_P(a, b)|=\Gamma(a, b)>0\]
for $(a, b)\in \partial\mathcal{H}^+_{\mathbb{R}}\cap U_{\delta}$, this completes the proof of Theorem \ref{thm:nondegenerate}. \QED

\medskip

\noindent
\textit{Proof of Corollary \ref{cor:ergodic}.}
The condition $\frac{1}{2}\log|\lambda_Q(a, b)| \ne \Lambda(a, b)$ in this paper is equivalent to the non-degeneracy assumption $\frac{1}{2}\lambda^u(\delta_Q) \ne \lambda^u_m$ in Theorem A (a) of~\cite{T}. Therefore, Corollary \ref{cor:ergodic} in Section 1 follows. \QED

\begin{figure}[th]
\includegraphics[width=9cm]{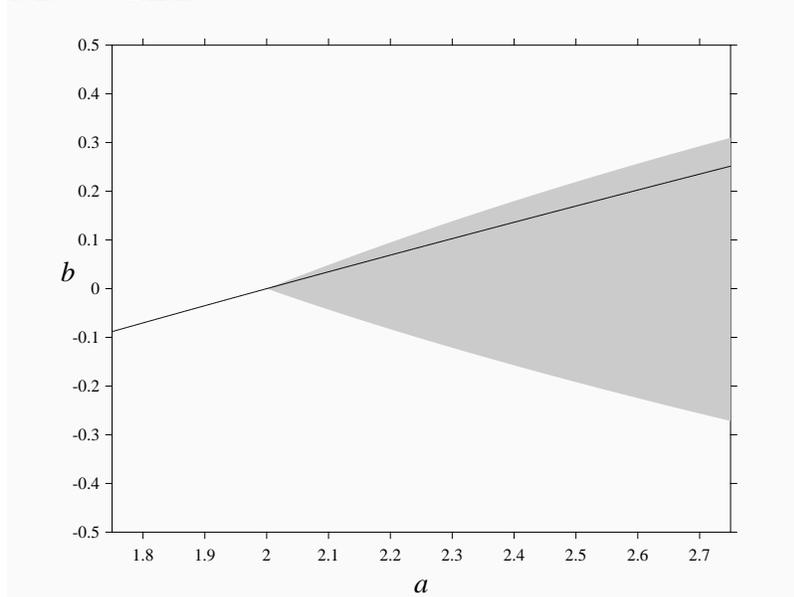}
\caption{Zero locus $\big\{(a, b) : \Gamma(a, b)=0\big\}$ for the fixed point $P$ (solid curve) and the hyperbolic horseshoe locus $\mathcal{H}_{\mathbb{R}}$ (shaded region).}
\label{FIG:per1}
\end{figure}

\begin{figure}[th]
\includegraphics[width=9cm]{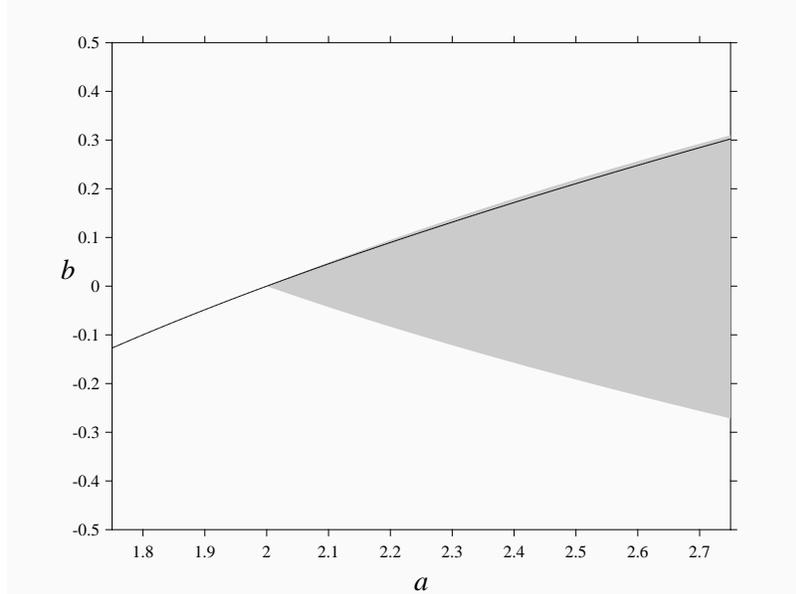}
\caption{Zero locus $\big\{(a, b) : \Gamma(a, b)=0\big\}$ for the period two cycle (solid curve) and the hyperbolic horseshoe locus $\mathcal{H}_{\mathbb{R}}$ (shaded region).}
\label{FIG:per2}
\end{figure}

\begin{figure}[th]
\includegraphics[width=9cm]{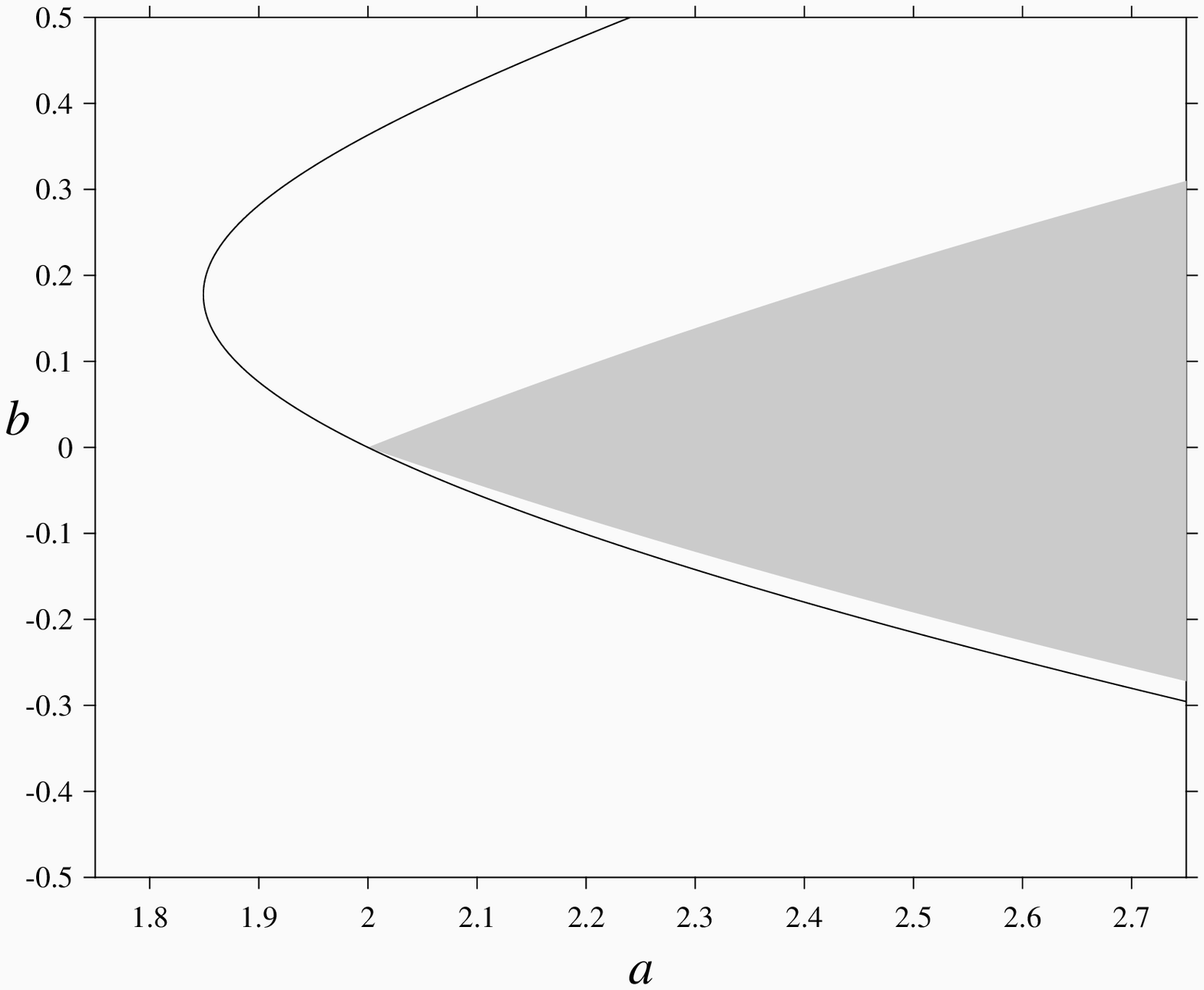}
\includegraphics[width=9cm]{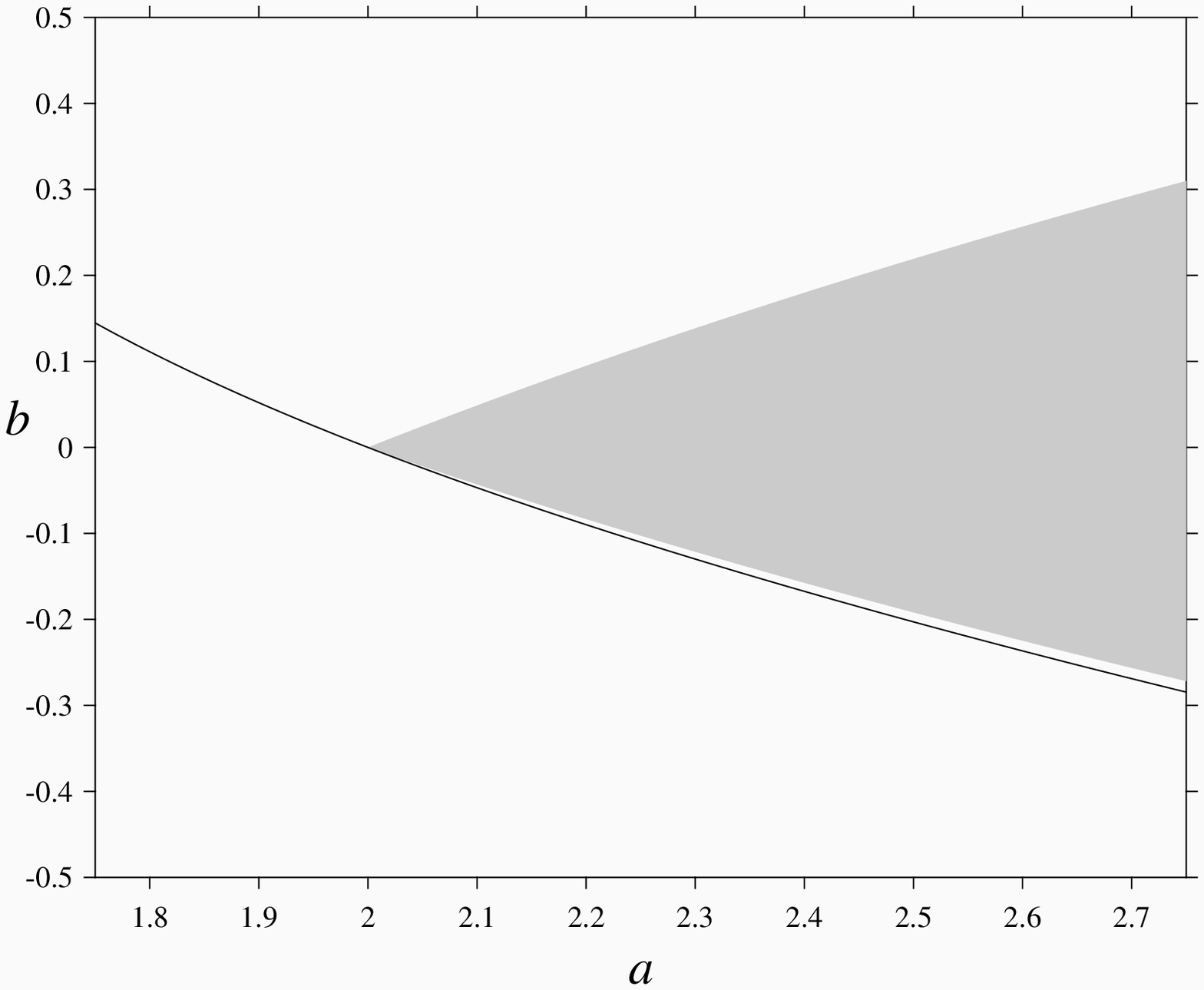}
\caption{Zero loci $\big\{(a, b) : \Gamma(a, b)=0\big\}$ for period three cycles (solid curves) and the hyperbolic horseshoe locus $\mathcal{H}_{\mathbb{R}}$ (shaded regions).}
\label{FIG:per3}
\end{figure}

In the proof above we only treated the fixed point $P(a, b)$ to estimate $\Lambda(a, b)$. If we could directly prove $\frac{1}{2}\log|\lambda_Q(2, 0)|>\Lambda_{\mu}(f_{2, 0})$ for the invariant measure $\mu$ supported on certain periodic orbit of $f_{2, 0}$ rather than $P(2, 0)$, this would imply $\frac{1}{2}\log|\lambda_Q(a, b)| \ne \Lambda(a, b)$ without using Theorem \ref{thm:slope}. However, this is not the case. To see this, let us recall that the Chebyshev map $p(x)=x^2-2$ on $[-2, 2]$ is topologically conjugate to the tent map $T(x)=2|x|-2$ on $[-2, 2]$ and the conjugacy map is differentiable except at the endpoints of $[-2, 2]$. This yields $\Lambda_{\mu}(f_{2, 0})=\log 2$ for the invariant measure $\mu$ supported on a periodic orbit except for the fixed point $Q(2, 0)$, and therefore $\frac{1}{2}\log|\lambda_Q(2, 0)|=\Lambda_{\mu}(f_{2, 0})$ always holds. This suggests that, in a sense, Theorem \ref{thm:slope} is necessary to prove Theorem \ref{thm:nondegenerate}.

For any periodic cycle of $f_{a, b}$, one can define $\Gamma(a, b)$ by replacing $\lambda_P(a, b)$ in its formula with the unstable eigenvalue of the cycle. Figure \ref{FIG:per1} is a picture where the cycle is chosen to be the fixed point $P(a, b)$. The upper-left part to the solid curve is where $\Gamma(a, b)>0$. We see that $\Gamma(a, b)>0$ for $(a, b)\in \partial\mathcal{H}^+_{\mathbb{R}}\cap U_{\delta}$ as was proved in Theorem \ref{thm:nondegenerate}. Figure \ref{FIG:per2} is a picture for the cycle of period two. The upper-left part to the solid curve is where $\Gamma(a, b)>0$. We see that the zero locus of $\Gamma(a, b)$ is almost identical to $\partial\mathcal{H}^+_{\mathbb{R}}\cap U_{\delta}$, hence it seems extremely hard to prove that $\Gamma(a, b)>0$ for $(a, b)\in \partial\mathcal{H}^+_{\mathbb{R}}\cap U_{\delta}$ by our method. Figure \ref{FIG:per3} is for the two cycles of period three. The upper-right part to the solid curve is where $\Gamma(a, b)<0$. We see that $\Gamma(a, b)<0$ for $(a, b)\in \partial\mathcal{H}^+_{\mathbb{R}}\cap U_{\delta}$. Hence, these numerical results suggest that the fixed point $P(a, b)$ is the only appropriate choice to prove Theorem \ref{thm:nondegenerate} among periodic cycles of period up to three.

\begin{table}
\[a^+_{50}=5.699311, \quad s^+_{50}=5.699311; \qquad \qquad a^+_{25}=3.371399, \quad s^+_{25}=3.632066;\]
\[a^+_{49}=5.586166, \quad s^+_{49}=5.615209; \qquad \qquad a^+_{24}=3.299545, \quad s^+_{24}=3.553407;\]
\[a^+_{48}=5.474702, \quad s^+_{48}=5.531115; \qquad \qquad a^+_{23}=3.229258, \quad s^+_{23}=3.475373;\]
\[a^+_{47}=5.364921, \quad s^+_{47}=5.447037; \qquad \qquad a^+_{22}=3.160526, \quad s^+_{22}=3.398010;\]
\[a^+_{46}=5.256821, \quad s^+_{46}=5.362984; \qquad \qquad a^+_{21}=3.093333, \quad s^+_{21}=3.321375;\]
\[a^+_{45}=5.150401, \quad s^+_{45}=5.278965; \qquad \qquad a^+_{20}=3.027665, \quad s^+_{20}=3.245521;\]
\[a^+_{44}=5.045662, \quad s^+_{44}=5.194992; \qquad \qquad a^+_{19}=2.963507, \quad s^+_{19}=3.170504;\]
\[a^+_{43}=4.942601, \quad s^+_{43}=5.111077; \qquad \qquad a^+_{18}=2.900839, \quad s^+_{18}=3.096385;\]
\[a^+_{42}=4.841218, \quad s^+_{42}=5.027232; \qquad \qquad a^+_{17}=2.839645, \quad s^+_{17}=3.023223;\]
\[a^+_{41}=4.741511, \quad s^+_{41}=4.943472; \qquad \qquad a^+_{16}=2.779903, \quad s^+_{16}=2.951083;\]
\[a^+_{40}=4.643479, \quad s^+_{40}=4.859812; \qquad \qquad a^+_{15}=2.721594, \quad s^+_{15}=2.880032;\]
\[a^+_{39}=4.547118, \quad s^+_{39}=4.776270; \qquad \qquad a^+_{14}=2.664695, \quad s^+_{14}=2.810135;\]
\[a^+_{38}=4.452427, \quad s^+_{38}=4.692862; \qquad \qquad a^+_{13}=2.609181, \quad s^+_{13}=2.741464;\]
\[a^+_{37}=4.359403, \quad s^+_{37}=4.609608; \qquad \qquad a^+_{12}=2.555027, \quad s^+_{12}=2.674089;\]
\[a^+_{36}=4.268042, \quad s^+_{36}=4.526530; \qquad \qquad a^+_{11}=2.502208, \quad s^+_{11}=2.608086;\]
\[a^+_{35}=4.178340, \quad s^+_{35}=4.443650; \qquad \qquad a^+_{10}=2.450694, \quad s^+_{10}=2.543529;\]
\[a^+_{34}=4.090294, \quad s^+_{34}=4.360992; \qquad \qquad a^+_{09}=2.400457, \quad s^+_{09}=2.480498;\]
\[a^+_{33}=4.003899, \quad s^+_{33}=4.278584; \qquad \qquad a^+_{08}=2.351464, \quad s^+_{08}=2.419074;\]
\[a^+_{32}=3.919149, \quad s^+_{32}=4.196452; \qquad \qquad a^+_{07}=2.303682, \quad s^+_{07}=2.359336;\]
\[a^+_{31}=3.836039, \quad s^+_{31}=4.114629; \qquad \qquad a^+_{06}=2.257078, \quad s^+_{06}=2.301374;\]
\[a^+_{30}=3.754561, \quad s^+_{30}=4.033143; \qquad \qquad a^+_{05}=2.211615, \quad s^+_{05}=2.245277;\]
\[a^+_{29}=3.674710, \quad s^+_{29}=3.952030; \qquad \qquad a^+_{04}=2.167254, \quad s^+_{04}=2.191140;\]
\[a^+_{28}=3.596478, \quad s^+_{28}=3.871326; \qquad \qquad a^+_{03}=2.123956, \quad s^+_{03}=2.139062;\]
\[a^+_{27}=3.519854, \quad s^+_{27}=3.791071; \qquad \qquad a^+_{02}=2.081677, \quad s^+_{02}=2.089152;\]
\[a^+_{26}=3.444831, \quad s^+_{26}=3.711303; \qquad \qquad a^+_{01}=2.040374, \quad s^+_{01}=2.041528;\]
\[a^+_{00}=2.000000, \quad s^+_{00}=2.000001.\]

\vspace{0.3cm}
\caption{Data of the values $a^+_k$ and $s^+_k$ for $0\leq k\leq 50$.}
\label{tab:aprx1}
\end{table}

\begin{table}
\[a^-_{50}=6.198261, \quad s^-_{50}=-6.198261; \qquad \qquad a^-_{25}=3.602401, \quad s^-_{25}=-4.187841;\]
\[a^-_{49}=6.075102, \quad s^-_{49}=-6.117656; \qquad \qquad a^-_{24}=3.519443, \quad s^-_{24}=-4.107984;\]
\[a^-_{48}=5.953555, \quad s^-_{48}=-6.037052; \qquad \qquad a^-_{23}=3.438081, \quad s^-_{23}=-4.028218;\]
\[a^-_{47}=5.833620, \quad s^-_{47}=-5.956450; \qquad \qquad a^-_{22}=3.358313, \quad s^-_{22}=-3.948552;\]
\[a^-_{46}=5.715297, \quad s^-_{46}=-5.875852; \qquad \qquad a^-_{21}=3.280138, \quad s^-_{21}=-3.868996;\]
\[a^-_{45}=5.598586, \quad s^-_{45}=-5.795257; \qquad \qquad a^-_{20}=3.203553, \quad s^-_{20}=-3.789558;\]
\[a^-_{44}=5.483487, \quad s^-_{44}=-5.714668; \qquad \qquad a^-_{19}=3.128555, \quad s^-_{19}=-3.710249;\]
\[a^-_{43}=5.369999, \quad s^-_{43}=-5.634087; \qquad \qquad a^-_{18}=3.055142, \quad s^-_{18}=-3.631080;\]
\[a^-_{42}=5.258123, \quad s^-_{42}=-5.553515; \qquad \qquad a^-_{17}=2.983311, \quad s^-_{17}=-3.552064;\]
\[a^-_{41}=5.147859, \quad s^-_{41}=-5.472953; \qquad \qquad a^-_{16}=2.913058, \quad s^-_{16}=-3.473213;\]
\[a^-_{40}=5.039205, \quad s^-_{40}=-5.392406; \qquad \qquad a^-_{15}=2.844381, \quad s^-_{15}=-3.394538;\]
\[a^-_{39}=4.932162, \quad s^-_{39}=-5.311873; \qquad \qquad a^-_{14}=2.777275, \quad s^-_{14}=-3.316056;\]
\[a^-_{38}=4.826730, \quad s^-_{38}=-5.231358; \qquad \qquad a^-_{13}=2.711737, \quad s^-_{13}=-3.237781;\]
\[a^-_{37}=4.722908, \quad s^-_{37}=-5.150863; \qquad \qquad a^-_{12}=2.647763, \quad s^-_{12}=-3.159728;\]
\[a^-_{36}=4.620695, \quad s^-_{36}=-5.070390; \qquad \qquad a^-_{11}=2.585347, \quad s^-_{11}=-3.081910;\]
\[a^-_{35}=4.520092, \quad s^-_{35}=-4.989944; \qquad \qquad a^-_{10}=2.524484, \quad s^-_{10}=-3.004349;\]
\[a^-_{34}=4.421097, \quad s^-_{34}=-4.909529; \qquad \qquad a^-_{09}=2.465171, \quad s^-_{09}=-2.927058;\]
\[a^-_{33}=4.323711, \quad s^-_{33}=-4.829147; \qquad \qquad a^-_{08}=2.407400, \quad s^-_{08}=-2.850056;\]
\[a^-_{32}=4.227931, \quad s^-_{32}=-4.748801; \qquad \qquad a^-_{07}=2.351167, \quad s^-_{07}=-2.773361;\]
\[a^-_{31}=4.133758, \quad s^-_{31}=-4.668497; \qquad \qquad a^-_{06}=2.296464, \quad s^-_{06}=-2.696986;\]
\[a^-_{30}=4.041191, \quad s^-_{30}=-4.588241; \qquad \qquad a^-_{05}=2.243285, \quad s^-_{05}=-2.620952;\]
\[a^-_{29}=3.950228, \quad s^-_{29}=-4.508035; \qquad \qquad a^-_{04}=2.191623, \quad s^-_{04}=-2.545274;\]
\[a^-_{28}=3.860869, \quad s^-_{28}=-4.427885; \qquad \qquad a^-_{03}=2.141471, \quad s^-_{03}=-2.469965;\]
\[a^-_{27}=3.773113, \quad s^-_{27}=-4.347800; \qquad \qquad a^-_{02}=2.092822, \quad s^-_{02}=-2.395034;\]
\[a^-_{26}=3.686957, \quad s^-_{26}=-4.267782; \qquad \qquad a^-_{01}=2.045667, \quad s^-_{01}=-2.320487;\]
\[a^-_{00}=2.000000, \quad s^-_{00}=-2.246320.\]

\vspace{0.3cm}
\caption{Data of the values $a^-_k$ and $s^-_k$ for $0\leq k\leq 50$.}
\label{tab:aprx2}
\end{table}


\begin{thebibliography}{[BS]}

\bibitem[AI]{AI} 
Z.~Arai, Y.~Ishii,
\textit{On parameter loci of the H\'enon family.}
Preprint available at \texttt{arXiv:1501.01368} (2015).

\bibitem[AM]{AM} 
Z.~Arai, K.~Mischaikow, 
\textit{Rigorous computations of homoclinic tangencies.}
SIAM Journal on Applied Dynamical Systems \textbf{5} (2006), 280--292.

\bibitem[BS]{BS} 
E.~Bedford, J.~Smillie, 
\textit{Real polynomial diffeomorphisms with maximal entropy: II. Small Jacobian.} 
Ergodic Theory Dynam. Systems \textbf{26} (2006), no. 5, 1259--1283.

\bibitem[DN]{DN} R.~Devaney, Z.~Nitecki,
\textit{Shift automorphisms in the H\'enon mapping.}
Comm. Math. Phys. \textbf{67} (1979), no. 2, 137--146. 

\bibitem[FM]{FM} 
S.~Friedland, J.~Milnor,
\textit{Dynamical properties of plane polynomial automorphisms.}
Ergodic Theory Dynam. Systems \textbf{9} (1989), no. 1, 67--99. 

\bibitem[L]{L}
R.~Leplaideur,
\textit{Thermodynamic formalism for a family of non-uniformly hyperbolic horseshoes and the unstable Jacobian.}
Ergodic Theory Dynam. Systems \textbf{31} (2011), no. 2, 423--447.

\bibitem[M]{M} J.~Milnor,
\textit{Nonexpansive H\'enon maps.}
Adv. Math. \textbf{69} (1988),  no. 1, 109--114.

\bibitem[ST]{ST} 
S.~Senti and H.~Takahasi, 
\textit{Equilibrium measures for the H\'enon map at the first bifurcation.} 
Nonlinearity \textbf{26} (2013), 1719--1741.

\bibitem[T]{T} 
H.~Takahasi, 
\textit{Equilibrium measures at temperature zero for H\'enon-like maps at the first bifurcation.} 
SIAM Journal on Applied Dynamical Systems \textbf{15} (2016), 106--124.

\end{thebibliography}
\end{document}